\documentclass[11pt,reqno]{amsart}
\usepackage{amssymb}
\usepackage[T1]{fontenc}
\usepackage{lmodern}
\usepackage{booktabs}
\usepackage{microtype}
\usepackage{tikz}
\usetikzlibrary{arrows.meta}
\usepackage[colorlinks=true,allcolors=blue!55!black]{hyperref}
\newtheorem{theorem}{Theorem}[section]
\newtheorem{lemma}[theorem]{Lemma}

\newcommand{\Sph}{\mathbb S}
\newcommand{\F}{\mathbb F_3}
\newcommand{\TMF}{\mathrm{TMF}}
\newcommand{\Jper}{\mathrm J^2}
\newcommand{\Jconn}{\mathrm j^2}
\newcommand{\BP}{\mathrm{BP}}
\newcommand{\Ext}{\operatorname{Ext}}
\newcommand{\gr}{\operatorname{gr}}
\newcommand{\unitrel}{\mathrel{\doteq}}
\DeclareRobustCommand{\xclass}{x_{153,3}}

\makeatletter
\@ifundefined{subjclassname@2020}{%
  \@namedef{subjclassname@2020}{\textup{2020} Mathematics Subject Classification}%
}{}
\makeatother

\title[Permanence and periodicity in the 153-stem]
{Permanence and \texorpdfstring{$144$}{144}-periodicity of\\
\texorpdfstring{$\xclass$}{x153,3} in the stable stems at the prime three}
\author{Zheyan Wan}
\address{Beijing Institute of Mathematical Sciences and Applications,
  Beijing 101408, China}
\email{wanzheyan@bimsa.cn}

\begin{document}

\begin{abstract}
The class $\xclass$ lies in stable degree $153$ and filtration $3$
on the algebraic starting page of the $3$-primary Adams--Novikov
spectral sequence of the sphere. We show that it survives and begins
a family with period $144$, answering both parts of a question of
Carrick and Davies. A comparison with periodic topological modular
forms and a bound on the possible differential targets prove
permanence. A divisibility calculation for the mod-$3$ Moore
spectrum supplies a lift to a finite spectrum with a degree-$144$
self-map. The resulting iterates have nonzero image in the
height-two detector of Carrick and Davies.
\end{abstract}

\maketitle

\section{Introduction}

Carrick and Davies \cite[Question~3.15]{CD} ask whether $\xclass$,
a class in stable degree $153$ and filtration $3$ on the initial page
of the $3$-primary Adams--Novikov spectral sequence of the sphere
\cite[\S1]{Nov}, survives to stable homotopy and begins a
$144$-periodic family. We prove that the class survives and that each
of its homotopy representatives begins a family of order-three
classes in degrees $153+144k$ for integers $k\geq0$. We construct the
family by lifting the representative to a finite spectrum and
iterating a degree-$144$ $v_2^9$ self-map of Belmont and Shimomura
\cite[Corollary~4.7]{BS}.

Write $\Sph$ for the sphere spectrum and $\Sph_{(3)}$ for its
localization at the prime $3$, which inverts integers prime to $3$.
For a spectrum $Y$, $\Sigma^n Y$ denotes its $n$-fold suspension
(desuspension for negative $n$) and $\pi_nY=[\Sigma^n\Sph,Y]$
its \(n\)th stable homotopy group.
Here $[A,B]$ denotes stable homotopy classes of maps from $A$ to
$B$, and $\pi_*Y$ denotes the graded collection of homotopy groups.
An integer $n$ indexing \(\pi_n\) is called a \emph{stem}.
We work at the prime $3$ throughout.

Chromatic homotopy theory organizes periodic families in the stable
stems by height \cite[Chapters~1--2]{RavNil}. At height one, the image
of the stable $J$-homomorphism, which is induced by one-point
compactification of orthogonal transformations, is closely related
to the Adams self-map of the mod-$3$ Moore spectrum. For a positive
integer $m$, write $Y/m$ for the cofiber of multiplication by $m$ on a spectrum $Y$,
which fits into the sequence $Y\xrightarrow{m}Y\to Y/m\to\Sigma Y$.
For $m=3$ the last two maps are denoted by $q_0$ and
$\partial_0$, respectively.
The mod-$3$ Moore spectrum is $M=\Sph/3$. The Adams self-map
$v_1\colon\Sigma^4M\to M$ induces multiplication by the degree-$4$
generator of $K(1)$-homology, where $K(i)$ denotes the \(i\)th
Morava \(K\)-theory \cite[Example~2.4.1(ii)]{RavNil}. Adams constructs
this map using complex $K$-theory \cite[\S12, proof of Theorem~1.7]{AdamsJ}.
For an integer $j\geq1$, write
$M_j=\Sph/(3,v_1^j)$ for the generalized Moore spectrum given by the cofiber of
$v_1^j\colon\Sigma^{4j}M\to M$. A power of a self-map here uses
composition with the appropriate suspensions. The Smith--Toda
complex $V(1)$ is $M_1$ \cite[Corollary~4.3 and the example following it]{TodaV}.

The nilpotence theorem of Devinatz, Hopkins, and Smith detects
nilpotence of self-maps of finite spectra using complex bordism
\cite[Theorem~1]{DHS}. The periodicity theorem of Hopkins and Smith gives a
positive-degree $v_n$ self-map on each finite $3$-local spectrum
of type $n\geq1$ \cite[Theorem~9]{HS}. Here type $n$ means that
its Morava $K$-homology vanishes at heights below $n$ and is
nonzero at height $n$. A $v_n$ self-map induces an isomorphism in
$K(n)$-homology and a nilpotent map in every other Morava
$K$-homology theory.

At height two, let $v_2$ denote the degree-$16$ periodicity
generator of the coefficients of $K(2)$. A $v_2^9$ self-map
has degree $144$ and induces multiplication by $v_2^9$ on
$K(2)$-homology. Behrens and Pemmaraju constructed such a map on
$V(1)$ \cite{BP}. Belmont and Shimomura constructed such maps on
$M_j$ for $2\leq j\leq8$ and used them to obtain families in the
stable stems \cite[Corollary 4.7 and \S5]{BS}.

We use the Brown--Peterson form of the Adams--Novikov spectral
sequence, writing $\BP_*Y$ for the homology of a spectrum $Y$ with
respect to the Brown--Peterson spectrum \cite{BrownPeterson}.
The initial page is the comodule cohomology group
\begin{equation}\label{eq:anss-e2}
E_2^{f,t}(Y)=\Ext^{f,t}_{\BP_*\BP}(\BP_*,\BP_*Y),
\end{equation}
where $\BP_*=\BP_*\Sph$ is the coefficient ring and
$\BP_*\BP$ is the ring of cooperations \cite[\S4.4]{Rav}. The
cohomological degree $f$ is the Adams--Novikov filtration and
$t$ is the internal degree. We write
$\mathcal E_r(n,f;Y)=E_r^{f,n+f}(Y)$ for the group on page $E_r$
in stem $n=t-f$ and filtration $f$. An omitted spectrum means
the sphere. The differential $d_r$ on page $E_r$ sends
$(n,f)$ to $(n-1,f+r)$. A $d_r$-cycle is a class killed by $d_r$,
and a $d_r$-boundary is a class in its image. A permanent cycle
survives nontrivially to the limiting page $E_\infty$.
The spectral sequence of the sphere has a product satisfying the
Leibniz rule, and the spectral sequences of the Moore spectra are
modules over it \cite{Moss}\cite[Theorem~2.3.3]{Rav}.

Let $F^f\pi_nY$ be the subgroup of Adams--Novikov filtration at
least $f$. The limiting page gives the associated graded group
$\gr^f_{\mathrm{AN}}\pi_nY=F^f\pi_nY/F^{f+1}\pi_nY$ of the
$3$-local homotopy groups for connective $Y$
\cite[Theorems~2.2.3(b) and~2.2.13]{Rav}. Here $\mathrm{AN}$
stands for Adams--Novikov, and we take $Y$ to be $3$-local in
this notation. The sphere has finite homotopy groups in positive
stems, where its $3$-local and $3$-completed groups agree.
A class in $\pi_nY$ that maps
to a given class of $E_\infty$ is called a representative of
that class. A hidden extension by three occurs when a class
is killed by $3$ on $E_\infty$ but multiplying a representative
by $3$ gives a nonzero class in higher filtration.

The notation $v_i$ also denotes the usual chromatic generators
of $\BP_*$, of internal degree $2(3^i-1)$ for $i\geq1$.
For the algebraic Greek-letter constructions, see
\cite[Theorems~2.2 and~2.6]{MRW}.
The Greek-letter class $\alpha_1$ is the algebraic connecting
class of $v_1$ modulo $3$, in bidegree $(3,1)$.
The divided class $\beta_{i/j}$ is defined, where the construction
applies, by taking the connecting classes of $v_2^i$ modulo
$(3,v_1^j)$ for multiplication by $v_1^j$ and then by $3$.
Its bidegree is $(16i-4j-2,2)$, and we abbreviate $\beta_{i/1}$
to $\beta_i$ \cite[\S5]{BS}. Juxtaposition and superscripts on
classes denote multiplication and its iterated powers. We use the
same symbols for homotopy representatives when the classes are
permanent and we are discussing products in homotopy.

Write $\mathrm{tmf}$ for connective topological modular forms
\cite[\S4]{HopkinsTMF}. Davies constructed stable Adams operations
on topological modular forms and introduced a connective height-two
image-of-$J$ spectrum \cite[Theorem~B and \S3.3]{DavOps}. We use
the $3$-complete version $\Jconn$ of Carrick and Davies
\cite[Definition~2.2 and \S2.1]{CD}. In high degrees its homotopy
groups are calculated from the fiber of $\psi^2-1$ on $\mathrm{tmf}$,
where $\psi^2$ is the Adams operation indexed by $2$ and $1$ is the
identity map. The unit map
$\Sph\to\Jconn$ gives the Hurewicz homomorphism on stable homotopy groups,
whose image is called the Hurewicz image. A class is detected by $\Jconn$
when its image under this homomorphism is nonzero. Let $\TMF$ denote
$3$-local periodic topological modular forms, with $\psi^2$ as constructed
in \cite[Definition~2.1]{DavHecke}. Write $\Jper$ for the
homotopy equalizer of $\psi^2$ and the identity on $\TMF$, whose
underlying spectrum is the fiber of $\psi^2-1$
\cite[Definition~3.5]{Dav}.
We write $\partial$ for the connecting map
$\pi_{n+1}\TMF\to\pi_n\Jper$ from this fiber sequence.
We use the same symbol for its connective counterpart
$\pi_{n+1}\mathrm{tmf}\to\pi_n\Jconn$ in the high degrees where
$\Jconn$ agrees with the fiber.

In topological modular forms, $\alpha$ and $\beta$ denote the
images of $\alpha_1$ and $\beta_1$, or their spectral sequence
classes when working on a page. The modular discriminant $\Delta$
has weight $12$ and bidegree $(24,0)$. In homotopy, powers of
$\Delta$ with exponent divisible by three denote powers of the
degree-$72$ class detected by $\Delta^3$
\cite[Proposition~4.6]{HopkinsTMF}\cite[\S\S5--6]{Bau}.
We use $\Jper$ to compare
differentials and $\Jconn$ to detect the periodic family. The
modified Adams--Novikov spectral sequence of $\Jper$ uses the
same stem and filtration convention as the ordinary one.

Carrick and Davies identify $\xclass$ on the $E_2$-page by the relation
$\alpha_1\xclass\unitrel\beta_5\beta_{6/3}$, where $\unitrel$
denotes equality up to multiplication by a unit at the prime $3$
\cite[Remark~3.14]{CD}.
They report that its counterpart in the $K(2)$-local sphere,
obtained by localizing with respect to $K(2)$-homology, is
permanent \cite[discussion following Question~3.15]{CD}. Naturality
then constrains a possible differential on $\xclass$ to have target
with zero image on the corresponding local page. Permanence in the
sphere requires controlling those possible targets there.

For the periodicity question, let
$\partial_j\colon M_j\to\Sigma^{4j+1}M$ be the connecting map of
the cofiber sequence defining $M_j$ for $j\geq1$. Together with the
mod-$3$ boundary $\partial_0\colon M\to\Sigma\Sph$, it gives
\begin{equation}\label{eq:moore-boundaries}
\pi_{155+4j}M_j\xrightarrow{\partial_j}
\pi_{154}M\xrightarrow{\partial_0}\pi_{153}\Sph .
\end{equation}
To construct a periodic family from a representative $X$ of $\xclass$,
we need a lift $z\in\pi_{155+4j}M_j$ through these two boundaries.
Such a lift requires $3X=0$ and a choice of $y\in\pi_{154}M$ with
$\partial_0y=X$ and $v_1^jy=0$. Given $z$ and a $v_2^9$ self-map
$V\colon\Sigma^{144}M_j\to M_j$, we write $V^kz$ for the class
obtained by applying the suspended \(k\)-fold iterate of $V$ to
$z$ for each integer $k\geq0$. To show $v_2^9$-periodicity, we must
show that all the double boundaries $\partial_0\partial_j(V^kz)$
are nonzero.

Write $\mathbb Z/m$ for the cyclic group of order $m$.
\begin{samepage}
\begin{theorem}\label{thm:main}
The class $\xclass$ is a permanent cycle in the Adams--Novikov
spectral sequence for $\Sph_{(3)}$. In particular,
$\gr^3_{\mathrm{AN}}\pi_{153}\Sph_{(3)}\cong\mathbb Z/3$,
generated by the class detected by $\xclass$.
Every representative of $\xclass$ has order three and begins
a $v_2^9$-periodic family of order-three, filtration-three classes. The class
in stem $153+144k$ is detected by
$\partial(\beta\Delta^{6(k+1)})$ in $\pi_*\Jconn$.
\end{theorem}
\end{samepage}

Sections~\ref{sec:data} and \ref{sec:proof} compute the possible
differential targets and rule them out using comparison with $\Jper$,
multiplication in the sphere, and a finite-page vanishing line,
which bounds the bidegrees that can still be nonzero on a fixed page.
Section~\ref{sec:three} resolves the hidden extension by three.
Section~\ref{sec:lift} uses $\beta_1$-divisibility in the Moore
spectrum to choose a lift with $v_1^3y=0$.
Section~\ref{sec:detect} then detects every iterate in $\Jconn$.
These stems are the exceptional degrees in Carrick and Davies's
description of the Hurewicz image
\cite[Theorem A]{CD}.

\section{The relevant \texorpdfstring{$E_2$}{E2}-groups and products}
\label{sec:data}

We use revision \texttt{f530898} of the data of Belmont and Wang
\cite{Bel}, generated with Belmont's prime-$3$ adaptation of
Wang's \texttt{MinimalResolution} program \cite{MR}. Wang describes
the underlying minimal-resolution algorithm in \cite{Wang}.
The data extend through internal degree $185$, which contains every
group looked up in the proof. The program computes
the algebraic Novikov spectral sequence and its multiplicative
information, which determine the indicated groups and products on
the Adams--Novikov $E_2$-page.

Write $\F$ for the field of three elements. The notation
$[f-i]$ in Belmont's generator names labels the basis element
with index $i$ in homological degree $f$ of the chosen resolution.
The hyphen separates the two indices. A prefix such as
$v_2^6$ in $v_2^6[8-2]$ records the coefficient monomial of
the leading term used to name the generator. In these names
$v_0$ denotes $3$ and $v_i$ for $i\geq1$ denotes the
Brown--Peterson generator of that index.

\begin{lemma}\label{lem:data}
The groups in the following table have the indicated $\F$-bases,
with $\varnothing$ denoting the empty basis of the zero group.
Here $b$ is the generator \texttt{v3\string^1[12-5]} in
$(152,12)$ and $w$ is the generator \texttt{[4-8]} in $(92,4)$.
\[
\begin{array}{cl}
\toprule
(n,f) & \text{Basis of }\mathcal E_2(n,f)\\
\midrule
(153,3) & \xclass\\
(152,8) & \beta_1^2\beta_4\beta_5\\
(152,12) & b\\
(152,16) & \beta_1^6w,\ \beta_1^7\beta_{6/3}\\
(152,20) & \varnothing\\
(91,9) & \beta_1^3\alpha_1\beta_4\\
(91,13) & \varnothing\\
(90,18) & \beta_1^9\\
\bottomrule
\end{array}
\]
The product $\alpha_1 b$ in $\mathcal E_2(155,13)$ is nonzero and
the map given by multiplication by $\alpha_1$,
\begin{equation}\label{eq:alpha1-surjection}
\alpha_1\cdot-\colon\mathcal E_2(153,7)\longrightarrow
                       \mathcal E_2(156,8)
\end{equation}
is surjective.
\end{lemma}

\begin{proof}
The following are the generator names in Belmont's data \cite{Bel},
in the order of the named elements in the table. Coefficients
of the single-generator products are suppressed.
\[
\begin{array}{ll@{\qquad}ll}
\toprule
\xclass & \texttt{v0\string^1v1\string^1v2\string^1v3\string^2[3-2]}
  & \beta_1^2\beta_4\beta_5 & \texttt{v2\string^6[8-2]}\\
\beta_1^6w & \texttt{[16-9]}
  & \beta_1^7\beta_{6/3} & \texttt{v2\string^3[16-3]}\\
\beta_1^3\alpha_1\beta_4 & \texttt{v2\string^3[9-0]}
  & \beta_1^9 & \texttt{[18-0]}\\
\bottomrule
\end{array}
\]
To reproduce the lookups, initialize \texttt{ANSSData()} from
\texttt{beta.py} at the root of the data repository \cite{Bel}. The method
\texttt{in\_deg(n,f)} lists the generators in stem $n$ and
filtration $f$, and \texttt{three} multiplies a class by $3$.
The $\beta$-products use the multiplication tables for the mod-$3$
Moore spectrum and the algebraic boundary, as implemented in
\texttt{beta.py}. For a Greek-letter monomial, apply the corresponding
multiplication methods to the unit \texttt{[0-0]}. For a product such as
$\beta_1^6w$, apply \texttt{beta1(w,6)} to the named generator $w$.
For the last two claims the \texttt{alpha1} table gives
\begin{align}
\alpha_1 b &=
  \texttt{v2\string^1[13-9]}+2\texttt{v2\string^3[13-5]},\label{eq:alpha1-b}\\
\alpha_1\texttt{[7-11]} &= \texttt{[8-10]},\label{eq:alpha1-product-711}\\
\alpha_1\texttt{v3\string^1[7-8]} &= \texttt{v3\string^1[8-8]},\label{eq:alpha1-product-78}\\
\alpha_1\texttt{v2\string^6[7-3]} &= 0.\label{eq:alpha1-product-73}
\end{align}
The three inputs in $(153,7)$ form a basis by \cite{Bel}.
The two nonzero images form a basis in $(156,8)$, so the
kernel of this map is spanned by \texttt{v2\string^6[7-3]}.
\end{proof}

\section{Permanence}\label{sec:proof}

The $\BP$ Hopf algebroid $(\BP_*,\BP_*\BP)$ is concentrated
in internal degrees divisible by four. Thus in stem $n$
only filtrations $f$ with $n+f\equiv0\pmod4$ can occur.
Since $d_r$ changes the internal degree by $r-1$, the possible
lengths are $5,9,13,\ldots$, and $\mathcal E_5=\mathcal E_2$.
A differential hitting $\xclass$ would have to start in
filtration $3-r<0$. Outgoing differentials have target
in stem $152$ and filtration $3+r$. By
\cite[Proposition 11.6 and Theorem 12.2]{BHS},
the sphere $E_{10}$-page vanishes at $f>14+(n+9)/23$.
At $n=152$ the bound is $f>21$, and $d_{21}$ would land in
filtration $24$. This excludes $r\geq21$, and Lemma~\ref{lem:data}
excludes $r=17$. Figure~\ref{fig:sphere-targets} shows the remaining
possibilities $r=5,9,13$.

\begin{figure}[tbp]
\centering
\begin{tikzpicture}[
  x=1cm,y=.185cm,font=\small,>=Latex,
  candidate/.style={->,semithick,dashed,blue!65!black,
    shorten <=3pt,shorten >=6pt},
  generator/.style={circle,fill=black,inner sep=1.65pt}
]
  \fill[black!5] (.55,21.5) rectangle (8.8,25.65);
  \draw[black!25] (.55,21.5) -- (8.8,21.5);
  \draw[->] (.55,1) -- (.55,26);
  \draw[->] (.55,1) -- (8.9,1) node[right] {stem \(n\)};
  \node[above] at (.55,26) {filtration \(f\)};
  \foreach \f in {3,8,12,16,20,24} {
    \draw (.48,\f) -- (.62,\f);
    \node[left] at (.48,\f) {\(\f\)};
  }
  \foreach \x/\n in {4.5/152,8.1/153} {
    \draw (\x,.65) -- (\x,1.35);
    \node[below] at (\x,.65) {\(\n\)};
  }

  \foreach \f/\r in {8/5,12/9,16/13} {
    \draw[candidate] (8.1,3) -- (4.5,\f);
    \node[anchor=west,blue!65!black] at (1.04,\f) {\(d_{\r}\)};
  }
  \foreach \f/\r in {20/17,24/21} {
    \node at (4.5,\f) {\(\times\)};
    \node[anchor=west] at (1.04,\f) {\(d_{\r}\)};
  }

  \node[generator] at (8.1,3) {};
  \node[right] at (8.22,3) {\(\xclass\)};
  \node[generator] at (4.5,8) {};
  \node[generator] at (4.5,12) {};
  \node[generator] at (4.37,16) {};
  \node[generator] at (4.63,16) {};
  \node[left] at (4.18,8) {\(\beta_1^2\beta_4\beta_5\)};
  \node[left] at (4.18,12) {\(b\)};
  \node[left] at (4.18,16) {\(\beta_1^6w,\ \beta_1^7\beta_{6/3}\)};
  \node[left] at (4.18,20) {\(E_2=0\)};
  \node[left] at (4.18,24) {\(E_{10}=0\)};
\end{tikzpicture}
\caption{Possible outgoing differentials from \(\xclass\) in the sphere
Adams--Novikov spectral sequence. The \(E_2\) entries come from
Belmont's data \cite{Bel}. Each filled dot denotes one generator,
and the differential labels are beside their targets. Dashed arrows
show the three candidates checked below.
Crosses mark the zero groups that exclude \(d_{17}\) and \(d_{21}\)
on the indicated pages. The shaded band starts above filtration \(21\),
so the \(E_{10}\) vanishing line also excludes every \(d_r\) with \(r>21\).}
\label{fig:sphere-targets}
\end{figure}

\subsection{The \texorpdfstring{$d_5$}{d5} comparison}

Davies's filtered model of $\Jper$
\cite[Definition 3.5 and Propositions 3.7--3.8]{Dav} gives the
comparison from the sphere and the maps of spectral sequences
induced by the equalizer boundary $\partial$ and by the projection
$p\colon\Jper\to\TMF$. On pages written in stem and filtration
coordinates, $\partial$ has degree $(-1,1)$, meaning it sends
$(n,f)$ to $(n-1,f+1)$. The modified spectral sequence is also
supported in internal degrees divisible by four, so its
$E_5$-page equals its $E_2$-page as well. Write $\tau$ for the
synthetic suspension comparison map
\cite[Definition~4.27]{Pst}\cite[Definition~9.11]{BHS}.
Multiplication by $\tau$ gives the filtration tower. Its $E_2$-page
is the bigraded homotopy of the synthetic model tensored with the
cofiber of $\tau$ \cite[\S2.2, Equation~(2.13)]{CD}.
Tensoring preserves the fiber sequence, so its exact sequence
identifies the image of $\partial$ in bidegree $(n,f)$
with the cokernel of $\psi^2-1$ on
$\mathcal E_2(n+1,f-1;\TMF)$. The Adams--Novikov spectral sequence
for periodic $\TMF$ agrees with its descent spectral sequence
\cite[Corollary~4.15]{CDSyn}. In positive filtration the
$\TMF$ page is spanned by the monomials $\beta^a\Delta^j$
and $\alpha\beta^a\Delta^j$ in bidegrees $(10a+24j,2a)$ and
$(3+10a+24j,2a+1)$, respectively \cite[Equation (3.2)]{Dav}.

In bidegree $(153,3)$ the $\TMF$ term is zero. The boundary from
bidegree $(154,2)$ of $\TMF$ therefore spans the $\Jper$ group
in $(153,3)$, which is generated by $\partial(\beta\Delta^6)$.
Indeed, $\beta\Delta^6$ is fixed by $\psi^2$ on $E_2$, and is
the generator of that bidegree. On $\TMF$,
$d_5(\Delta)\unitrel\alpha\beta^2$
\cite[Equation~(6.1)]{Bau}\cite[Proposition~7.2 and Corollary~E]{CDN}.
Also $d_5(\beta)=0$ \cite[Figure~A.1 and Corollary~E]{CDN}.
The Leibniz
rule gives $d_5(\beta\Delta^6)\unitrel6\alpha\beta^3\Delta^5=0$
because $3\alpha=0$. Hence the image of $\xclass$ is a $d_5$-cycle.

By \cite[Propositions 3.10 and 3.11]{Dav}, the classes $\beta_1$ and
$\beta_4$ project to $\beta$ and $\beta\Delta^2$, respectively,
and $\beta_5$ maps to $\partial(\alpha\Delta^3)$, up to units.
On the synthetic models, the projection $p$ gives the action
of the equalizer on the $\TMF$ term. Since $\psi^2p\simeq p$,
the boundary is linear over the equalizer, and this linearity
descends to $E_2$. Consequently
the image of the generator of the possible
$d_5$ target group is
\begin{equation}\label{eq:d5-target-image}
 \beta_1^2\beta_4\beta_5\longmapsto
              \partial(\alpha\beta^3\Delta^5)
 \quad\text{up to a unit.}
\end{equation}
Here $\alpha\beta^3\Delta^5$ generates the $\TMF$ term in
$(153,7)$. It has order three, and $\psi^2$ fixes it, since
$\psi^2$ fixes $\alpha$ and $\beta$ and acts on $\Delta$ by $2^{12}$
\cite[Corollary~2.12]{DavHecke}\cite[Remark~3.4]{Dav}.
The equalizer exact sequence therefore shows that its boundary
is nonzero. The comparison is injective on $\mathcal E_5(152,8)$,
so naturality forces $d_5(\xclass)=0$.

\subsection{The \texorpdfstring{$d_9$}{d9} obstruction}

Suppose $d_9(\xclass)\ne 0$. Its target is represented by a unit
multiple of $b$. The group $\mathcal E_2(152,12)$ is
one-dimensional by \cite{Bel}, so $b$ can survive to $E_9$ only if
every $d_5$ from $(153,7)$ is zero. By Lemma~\ref{lem:data},
each element in $(156,8)$ has the form $\alpha_1 A$ with
$A\in\mathcal E_5(153,7)$. Since $\alpha_1$ is permanent, the
Leibniz rule gives $d_5(\alpha_1 A)=0$. Thus the only possible
$d_5$ into $(155,13)$ is zero. The product $\alpha_1b$ is nonzero
on $E_5$ and is a $d_5$-cycle because $b$ is, so it is still
nonzero on $E_9$.

On $E_2$, $\alpha_1\xclass\unitrel\beta_5\beta_{6/3}$ by
\cite[Remark~3.14]{CD}. The two factors on the right are permanent
cycles \cite[Corollary~1.2]{BP}\cite[Theorem~5.1]{BS}, so their product
supports no differential by the Leibniz rule. But the same rule
gives $d_9(\alpha_1\xclass)=\pm\alpha_1d_9(\xclass)\ne0$ under
our assumption, a contradiction.

\subsection{The \texorpdfstring{$d_{13}$}{d13} target group}

Ravenel's differential is
$d_9(\alpha_1\beta_4)\unitrel\beta_1^6$
\cite[Theorem 7.5.3]{Rav} (see also \cite[Lemma 12.10]{BHS}).
We check that $w$ is a $d_9$-cycle.
The class $\beta_1^9$ is a nonzero $d_5$-cycle, and its possible
$d_5$ source $\mathcal E_5(91,13)$ is zero by the $E_2$ data
\cite{Bel}. Hence it is nonzero on $E_9$. The Leibniz rule then gives
\begin{equation}\label{eq:d9-product}
d_9(\beta_1^3\alpha_1\beta_4)\unitrel\beta_1^9.
\end{equation}
In particular, the generator of $\mathcal E_5(91,9)$ survives to
$E_9$, so it cannot be a $d_5$-boundary. This is the only possible
$d_5$ target for $w$. The $d_9$ target for $w$ is itself zero by
Lemma~\ref{lem:data}, and filtration rules out incoming
differentials to $w$ through $E_9$.

The two basis elements $\beta_1^6w$ and $\beta_1^7\beta_{6/3}$
in $\mathcal E_2(152,16)$ are $d_5$-cycles, so their images span
$\mathcal E_9(152,16)$. On $E_9$ they are the $d_9$-boundaries
obtained by multiplying $d_9(\alpha_1\beta_4)$ by $w$ and by
$\beta_1\beta_{6/3}$, respectively. It follows that
$\mathcal E_{10}(152,16)=0$. This excludes $d_{13}(\xclass)$
and proves the permanence assertion of Theorem~\ref{thm:main}.

\section{The hidden extension by three}\label{sec:three}

Let $X\in\pi_{153}\Sph_{(3)}$ represent $\xclass$. In this
section we abbreviate $F^f\pi_{153}\Sph_{(3)}$ to $F^f$.
The relation $3\xclass=0$ on $E_2$ persists to $E_\infty$,
so $3X$ has zero image in $F^3/F^4$ and lies in $F^4$.
In stem $153$, sparseness requires $153+f\equiv0\pmod4$,
so $f\equiv3\pmod4$. The graded groups in filtrations $4$,
$5$, and $6$ therefore vanish, giving
\begin{equation}\label{eq:three-x-filtration}
  3X\in F^4=F^5=F^6=F^7 .
\end{equation}
Also $\alpha_1(3X)=0$ because $3\alpha_1=0$. The $E_{10}$
vanishing line in stem $153$ is $f>21+1/23$, leaving
$7$, $11$, $15$, and $19$ as candidates for the filtration
of a nonzero $3X$. The $E_2$ group in filtration $15$ is zero
by \cite{Bel}.

In filtration $7$, the $\alpha_1$-kernel on $E_5$ is spanned
by $z_7=\texttt{v2\string^6[7-3]}$, by the multiplication
calculation in Section~\ref{sec:data}. The product tables \cite{Bel}
identify it on $E_5$ with $\beta_1\beta_{6/3}(\alpha_1\beta_4)$
up to a unit. The factors are $d_5$-cycles, so the product is a
$d_5$-cycle. Its $d_9$ is $\beta_1^7\beta_{6/3}$ up to a unit
by the Leibniz rule. This differential is nonzero since
$\beta_1^7\beta_{6/3}$ maps to
$\beta^8\Delta^3$ in the Adams--Novikov spectral sequence of
$\TMF$ \cite[Proposition~3.10]{Dav}.
Belmont and Shimomura compute the corresponding
$K(2)$-local image of $\beta_{6/3}$ \cite[Proposition~6.4(3)]{BS}.
The monomial $\beta^8\Delta^3$ is a
$d_5$-cycle since $d_5\beta=0$ and $d_5(\Delta^3)=0$.
In filtration $11$ the only monomials are $\alpha\beta^5\Delta^j$,
in stems $53+24j$ by the formula in Section~\ref{sec:proof}.
Thus $\mathcal E_5(153,11;\TMF)=0$, and $\beta^8\Delta^3$ remains
nonzero on $E_9$.
In stem $156$ and filtration $8$, the only possible incoming
differential is $d_5$ from $(157,3)$, which is spanned by
$\alpha_1\beta_{10}$ by the $E_2$ data \cite{Bel}. This is a
permanent cycle because $\beta_{10}$ is permanent
\cite[Theorem~1]{OkaBeta}\cite[Corollary~1.2]{BP}.
Any nonzero $\alpha_1$-product of a
surviving class in $(153,7)$ thus survives in $(156,8)$.
Since $z_7$ supports $d_9$, multiplication by $\alpha_1$ is
injective on the filtration-$7$ associated graded.
The filtration-$7$ symbol of $3X$ is killed by $\alpha_1$ and
therefore vanishes. Thus $3X\in F^8=F^{11}$ by sparseness.

For filtration $11$, Belmont's data \cite{Bel} identify the
generator as $\beta_4 u$ up to a unit, with
$u=\texttt{[9-6]}\in\mathcal E_5(95,9)$. The
classical differential
$d_5(\beta_4)\unitrel\alpha_1\beta_1^2\beta_{3/3}$
\cite[\S6.4, before Lemma~6.4.12]{Rav}\cite[Table~1]{BS}
gives a nonzero component of
$d_5(\beta_4u)$ along $\beta_1^6w$. Indeed, the products in
\cite{Bel} show that $\alpha_1\beta_1^2\beta_{3/3}u$ is a unit multiple of
\texttt{[16-9]}. The only possible $d_5$ target for $u$ is
\texttt{[14-3]} in $(94,14)$, whose product with $\beta_4$ is a
unit multiple of the other generator \texttt{v2\string^3[16-3]}
in $(152,16)$. Thus in
$d_5(\beta_4u)=d_5(\beta_4)u+\beta_4d_5(u)$ the second term
cannot cancel the first, so the filtration-$11$ generator
supports $d_5$.

In filtration $19$, the data \cite{Bel} give the sole generator
$\texttt{[19-7]}\unitrel\beta_1^6 u'$, where
$u'=\texttt{[7-7]}$ lies in $(93,7)$. If $d_5u'$ is nonzero,
it is a unit multiple of \texttt{[12-4]} in $(92,12)$, and its
product with $\beta_1^6$ is the nonzero class
\texttt{[24-4]}. Hence $\beta_1^6u'$ supports $d_5$ in this case.
If $u'$ is a $d_5$-boundary, its product with $\beta_1^6$
is already zero on $E_9$. If $u'$ survives to $E_9$, then
its $d_9$ target in $(92,16)$ is zero. In this case
$d_9((\alpha_1\beta_4)u')\unitrel\beta_1^6u'$, so the
product is zero on $E_{10}$. This exhausts the possible
leading terms of $3X$. Strong convergence
\cite[Lemma 11.8]{BHS} gives a separated filtration, so $3X=0$.

\section{The generalized Moore lift}\label{sec:lift}

For $M=\Sph/3$, the mod-$3$ map $q_0\colon\Sph\to M$ includes
the bottom cell. The sphere element $\beta_1$ of degree $10$
acts on $M$ by a map $\Sigma^{10}M\to M$. In the composite
$\beta_1v_1^2$ we first apply the suspension of $v_1^2$ by
$10$ and then this action. Similarly, $\beta_1\pi_{148}M$
denotes the image of multiplication by $\beta_1$ on
$\pi_{148}M$.

\begin{lemma}\label{lem:moore}
The self-map $\beta_1 v_1^2\colon\Sigma^{18}M\to M$ is
nullhomotopic, and
\begin{equation}\label{eq:moore-divisibility}
F^6\pi_{158}M\ \subseteq\ \beta_1\pi_{148}M .
\end{equation}
\end{lemma}

\begin{proof}
Write $L_{K(1)}Y$ for localization of a spectrum $Y$ with
respect to $K(1)$-homology. The mod-$3$ cofiber sequence and
Toda's calculation \cite[Chapters~XIII--XIV]{TodaComp}
give $\pi_{18}M=0$ and $\pi_{19}M\cong\F$. The latter comes
from the image of $J$ \cite[Theorem~1.5]{AdamsJ} and maps
isomorphically to $\pi_{19}L_{K(1)}M$ by Ravenel's height-one
calculation \cite[Theorem~8.10(b)--(c) and Remark~8.12]{RavLoc},
which also gives $\pi_{18}L_{K(1)}M=0$.
For $N=M$ or $L_{K(1)}M$, mapping the
suspended Moore cofiber sequence into $N$ gives
\begin{equation}\label{eq:moore-mapping-sequence}
  \pi_{19}N\xrightarrow{3}\pi_{19}N
  \longrightarrow[\Sigma^{18}M,N]\longrightarrow\pi_{18}N=0.
\end{equation}
Thus $[\Sigma^{18}M,N]\cong\pi_{19}N/3$, and localization
induces an isomorphism between the two mapping groups.
The map $\beta_1v_1^2$ localizes to zero because
$\pi_{10}L_{K(1)}\Sph=0$.

For an abelian group $G$, write $G/3$ for its quotient by
$3G$ and $G[3]$ for its subgroup killed by $3$. The sequence
of $\BP_*\BP$-comodules for $M$ gives the short exact sequence
\begin{equation}\label{eq:mooreE2}
0\longrightarrow\mathcal E_2(n,f)/3
\xrightarrow{q_0}\mathcal E_2(n,f;M)
\xrightarrow{\partial_0}\mathcal E_2(n-1,f+1)[3]
\longrightarrow0 ,
\end{equation}
The maps in this sequence are induced by the bottom-cell inclusion
and the mod-$3$ boundary. They commute with $\beta_1$
multiplication with the indicated bidegree shift of $\partial_0$.
In Table~\ref{tab:moore}, the bottom group at $(n,f)$ is the left
term of \eqref{eq:mooreE2} and the top group is the right term,
coming from the top cell of $M$. Each generator is named in
Belmont's notation, and $0$ denotes the zero group.
We obtain the two end groups from \texttt{in\_deg} and
\texttt{three}, and use \texttt{beta1(g,k)} for multiplication of
a generator $g$ by $\beta_1^k$ in Belmont's data \cite{Bel}.
Reading across each pair of columns gives the product up to a unit. In particular
the first row is a surjection on the Moore $E_2$-page, and
rows two through six are isomorphisms.

\begin{table}[t]
\centering
\caption{The bottom and top groups and their $\beta_1$ products,
computed from Belmont's data \cite{Bel} using \eqref{eq:mooreE2}.
The target of a row with $\beta_1^k$ has bidegree $(n+10k,f+2k)$.}
\label{tab:moore}
\small
\setlength{\tabcolsep}{4pt}
\begin{tabular}{@{}ccllll@{}}
\toprule
 & Source & \multicolumn{2}{c}{Source generators}
 & \multicolumn{2}{c}{Target generators}\\
\cmidrule(lr){3-4}\cmidrule(l){5-6}
Action & $(n,f)$ & Bottom & Top & Bottom & Top\\
\midrule
$\beta_1$ & $(148,4)$ & $v_2^6[4-4]$ & $v_2^3v_3[5-3]$
 & $0$ & $v_2^3v_3[7-3]$\\
$\beta_1$ & $(148,8)$ & $v_3[8-7]$ & $0$ & $v_3[10-7]$ & $0$\\
$\beta_1$ & $(147,9)$ & $0$ & $v_2^6[10-0]$ & $0$ & $v_2^6[12-0]$\\
$\beta_1$ & $(147,13)$ & $v_3[13-4]$ & $v_3[14-3]$
 & $v_3[15-4]$ & $v_3[16-3]$\\
$\beta_1$ & $(147,17)$ & $0$ & $[18-7]$ & $0$ & $[20-7]$\\
$\beta_1^5$ & $(108,4)$ & $v_2^3[4-5]$ & $v_3[5-4]$
 & $v_2^3[14-5]$ & $v_3[15-4]$\\
\bottomrule
\end{tabular}
\end{table}

\begin{figure}[tbp]
\centering
\begin{minipage}[c]{.39\linewidth}
\centering
\small (a) Candidate differentials\par\smallskip
\begin{tikzpicture}[
  x=1cm,y=.27cm,font=\scriptsize,>=Latex,
  possible/.style={->,dashed,black!58,shorten <=4pt,shorten >=4pt},
  degree/.style={rectangle,fill=black,inner sep=1.65pt}
]
  \foreach \f in {4,8,9,13,17} {
    \draw[black!13] (.43,\f) -- (2.24,\f);
    \draw[black!13] (2.63,\f) -- (4.42,\f);
    \node[left] at (.38,\f) {\(\f\)};
  }
  \draw[->] (.43,2.53) -- (.43,18.05);
  \node[above,right] at (.43,17.95) {filtration \(f\)};
  \draw (.43,2.53) -- (2.24,2.53);
  \draw (2.63,2.53) -- (4.42,2.53);
  \node[below] at (2.42,.94) {stem \(n\)};
  \foreach \x/\n in {.94/147,1.91/148,3.13/147,4.1/148} {
    \draw (\x,2.3) -- (\x,2.76);
    \node[below] at (\x,2.3) {\(\n\)};
  }

  \foreach \f in {9,13,17} {
    \draw[possible] (1.91,4) -- (.94,\f);
    \node[degree] at (.94,\f) {};
  }
  \foreach \f in {13,17} {
    \draw[possible] (4.1,8) -- (3.13,\f);
    \node[degree] at (3.13,\f) {};
  }
  \node[degree] at (1.91,4) {};
  \node[degree] at (4.1,8) {};
\end{tikzpicture}
\end{minipage}\hfill
\begin{minipage}[c]{.57\linewidth}
\centering
\small (b) Multiplication by \(\beta_1\)\par\smallskip
\renewcommand{\arraystretch}{1.33}
\setlength{\tabcolsep}{5pt}
\begin{tabular}{@{}rcll@{}}
\toprule
Source & & Target & Page check\\
\midrule
\((148,4)\)  & {\color{blue!65!black}\(\longrightarrow\)}
  & \((158,6)\)  & surjective on \(E_2\)\\
\((148,8)\)  & {\color{blue!65!black}\(\longrightarrow\)}
  & \((158,10)\) & isomorphism on \(E_2\)\\
\((147,9)\)  & {\color{blue!65!black}\(\longrightarrow\)}
  & \((157,11)\) & injective on \(E_5\)\\
\((147,13)\) & {\color{blue!65!black}\(\longrightarrow\)}
  & \((157,15)\) & injective on \(E_9\)\\
\((147,17)\) & {\color{blue!65!black}\(\longrightarrow\)}
  & \((157,19)\) & injective on \(E_{13}\)\\
\bottomrule
\end{tabular}\par\smallskip
The \(d_5\) from \((158,14)\) is zero.
\end{minipage}
\caption{Moore spectrum bidegrees used in the
\(\beta_1\)-divisibility argument. The \(E_2\) groups and products
are computed from \cite{Bel} using \eqref{eq:mooreE2}. In (a), each
square marks a nonzero \(E_2\)-group in one of the two candidate differential
patterns. A dashed arrow from filtration \(f\) to filtration \(g\)
represents a possible \(d_{g-f}\). In (b), each arrow adds \(10\) to
the stem and \(2\) to the filtration. The page checks identify the
injectivity needed to compare differentials on the two sides.}
\label{fig:moore-pages}
\end{figure}

We first check that the source of the last row has permanent
representatives. The bottom-cell generator comes from the
permanent product $\beta_2\beta_{6/3}$ in $(108,4)$
\cite[Proposition~3.13]{CD}. Its image on $E_2$ is nonzero by
\eqref{eq:mooreE2} and Table~\ref{tab:moore}, and no differential
can hit filtration four.
For the top-cell generator, let $c$ be the sphere class
$v_3[5-4]$ in $(107,5)$. Ravenel's table \cite[Table~A3.4]{Rav}
gives two independent classes modulo the image of $J$ in stem
$107$, namely $\beta_2\gamma_2$ and $\beta_2x_{81}$. Here
$\gamma_2$ is the height-three Greek-letter element in stem $81$
and $x_{81}$ is the other generator in that stem in Ravenel's
table. The two products in stem $107$ are killed by $3$.
Together with the image of $J$ \cite[Theorem~1.5]{AdamsJ}
they give the direct sum decomposition
$\pi_{107}\Sph_{(3)}\cong\mathbb Z/81\oplus(\mathbb Z/3)^2$.
The cyclic order-$81$ filtration-one group survives by
\cite[Theorem 5.3.7(a)]{Rav}.
The $E_2$ data \cite{Bel} leave three other possible filtrations
below the $E_{10}$ vanishing line, namely $5$, $13$, and $17$,
each of dimension one.
In filtration $17$ the generator is
$\beta_1^6[5-3]$ up to a unit, for $[5-3]$ in $(47,5)$.
Its potential $d_5$ and $d_9$ targets $(46,10)$ and $(46,14)$
are zero in the same data. It follows that $\beta_1^6[5-3]$
vanishes on $E_{10}$, since $\beta_1^6$ is a $d_9$-boundary.
The known group order $81\cdot3^2$ then forces both
filtration-$5$ and filtration-$13$ classes to survive.
The filtration-one quotient restricts to an isomorphism on
the cyclic order-$81$ summand, since the order-three
summands have zero image modulo $3$ in that quotient.
Its kernel therefore has exponent three. In particular,
$c$ has an order-three representative that lifts to $\pi_{108}M$.

Multiplication by $3$ is injective on $\BP_*$, so the Moore
cofiber sequence gives
$0\to\BP_*\xrightarrow{3}\BP_*\to\BP_*M\to 0$ and
$\BP_*\partial_0=0$. The geometric boundary theorem
of Bruner \cite{BrunerBoundary}, in the form of
\cite[Theorem~2.3.4(a)--(c)]{Rav}, gives maps of bidegree
$(-1,1)$ on all pages from $E_2$ onward. The map on $E_2$
is the Ext connecting map in \eqref{eq:mooreE2}, and the map
on $E_\infty$ is the associated graded map of $\partial_0$
with filtration shift one. In particular $\partial_0$ raises
filtration by at least one. For the earlier permanent-cycle
formulation, see \cite[Theorem~1.7]{JMWZ}.
A lift of filtration zero would
have nonzero boundary in filtration one. Indeed, by
\eqref{eq:mooreE2} and \cite{Bel} the algebraic boundary is
injective in $(108,0)$, and its target is on the permanent
filtration-one line. A lift in $F^8\pi_{108}M$ would have
boundary in $F^9\pi_{107}\Sph_{(3)}$. Since the boundary
represents $c$ in filtration five, sparseness forces the
lift to have filtration four, with leading symbol mapping to $c$.
The two representatives span $\mathcal E_2(108,4;M)$, since
their symbols give the bottom generator and a lift of the top
generator in \eqref{eq:mooreE2}. Their products with $\beta_1^5$
span $\mathcal E_2(158,14;M)$ by the last row of
Table~\ref{tab:moore}. Consequently the $d_5$ from $(158,14)$
is zero, and every surviving class there is represented by such a
homotopy product.

To lift the first two rows to homotopy, we need the later-page
injections shown in Figure~\ref{fig:moore-pages}. For the
row with source $(147,13)$ the map is injective on $E_9$.
All possible $d_5$-boundaries in the target $(157,15)$
come from $(158,10)$, which on $E_5$ is in the image of
$\beta_1$ by the second row. Thus if a $d_5$-cycle maps to
a $d_5$-boundary, the Leibniz rule and injectivity on $E_5$
show it was already a boundary.

The row from $(147,17)$ to $(157,19)$ is injective on $E_9$,
since the only possible $d_5$-sources at the target are
in $(158,14)$, where we have shown $d_5=0$.
It remains injective on $E_{13}$. A possible $d_9$-boundary in
$(157,19)$ comes from a class in $(158,10)$. By the second row
of the table we can choose a preimage $B$ in $(148,8)$ on $E_5$
such that $\beta_1 B$ represents this class. Since
$d_5(\beta_1 B)=0$, injectivity on $E_5$ of the row from
$(147,13)$ gives $d_5B=0$. Hence the $d_9$-boundary in question
is $\beta_1 d_9 B$, and injectivity on $E_9$ shows that any
cycle mapping to this boundary was already a $d_9$-boundary.

A surviving class in $(158,10)$ now has a surviving antecedent
in $(148,8)$. The just-proved injections force the antecedent's
$d_5$ and $d_9$ to vanish, and an incoming differential would
kill its product with $\beta_1$ as well. For a survivor in
$(158,6)$, an antecedent in $(148,4)$ is a $d_5$, $d_9$, and
$d_{13}$ cycle by the corresponding injections from $(147,9)$,
$(147,13)$, and $(147,17)$. The vanishing line
\begin{equation}\label{eq:moore-vanishing-line}
  \mathcal E_{10}(n,f;M)=0 \qquad\text{for } f>14+(n+9)/23
\end{equation}
holds for $M$ because the line for the sphere is strong, so it
applies after tensoring with the $\BP$-synthetic lift of a connective
spectrum \cite[Definition 11.2(4), Proposition 11.6, and Theorem 12.2]{BHS}.
It excludes all remaining outgoing differentials from the antecedents,
whose next target filtration is at least $21$ in stem $147$.

In stem $158$, the two end groups of \eqref{eq:mooreE2}
are zero at filtration $18$ by \cite{Bel}. Filtrations at least
$22$ are above the $E_{10}$ line. By sparseness, the only possible
nonzero associated graded groups from filtration $6$ onward
occur at $6$, $10$, and $14$. Multiplication by $\beta_1$
surjects onto the first two from filtrations $4$ and $8$ in
stem $148$. In filtration $14$, the last row of
Table~\ref{tab:moore} represents each graded class by a product
$\beta_1^5W$ with $W\in F^4\pi_{108}M$. We can take its
$\beta_1$-antecedent to be $\beta_1^4W\in F^{12}\pi_{148}M$.
Thus for each $f\in\{6,10,14\}$ the map
\begin{equation}\label{eq:graded-beta1}
  \beta_1\colon \gr_{\mathrm{AN}}^{f-2}\pi_{148}M
      \longrightarrow \gr_{\mathrm{AN}}^f\pi_{158}M
\end{equation}
is surjective.

Let $A\in F^6\pi_{158}M$. Choose $U_6\in F^4\pi_{148}M$
such that $\beta_1U_6$ and $A$ agree in
$F^6\pi_{158}M/F^7\pi_{158}M$. The difference lies in
$F^7\pi_{158}M=F^{10}\pi_{158}M$. Use the same construction
in filtrations $10$ and $14$ to choose
$U_{10}\in F^8\pi_{148}M$ and $U_{14}\in F^{12}\pi_{148}M$.
Then $A-\beta_1(U_6+U_{10}+U_{14})$ lies in
$F^{15}\pi_{158}M=F^{22}\pi_{158}M$. All graded quotients
from filtration $22$ onward vanish, and strong convergence
\cite[Lemma 11.8]{BHS} gives $\bigcap_f F^f\pi_{158}M=0$.
Hence $F^{22}\pi_{158}M=0$ and
$A=\beta_1(U_6+U_{10}+U_{14})$, as required.
\end{proof}

By Section~\ref{sec:three}, we may choose $y\in\pi_{154}M$
with $\partial_0y=X$. The only possible filtration below
$6$ of $v_1y\in\pi_{158}M$ is $2$. By \eqref{eq:mooreE2}
and Belmont's data \cite{Bel}, the group $\mathcal E_2(158,2;M)$
has zero bottom and top generated by $\alpha_1\beta_{10}$
in $(157,3)$. Taking underlying spectra in
\cite[Lemma~2.6]{CD} gives
$\partial_0 v_1q_0(\beta_{10})=\alpha_1\beta_{10}$
up to a unit. The boundary has a nonzero filtration-three
symbol because both factors of $\alpha_1\beta_{10}$ are permanent
and no differential can hit this filtration. Since $\partial_0$
raises filtration by at least one, $v_1q_0(\beta_{10})$ has
filtration at most two. Sparseness forces exact filtration two.
Hence $v_1q_0(\beta_{10})$ represents a generator of the
filtration-two quotient. Subtracting a multiple of
$q_0(\beta_{10})$ from $y$ therefore gives a lift still
mapping to $X$ and satisfying $v_1y\in F^6\pi_{158}M$.
By Lemma~\ref{lem:moore} we can write $v_1y=\beta_1 U$
with $U\in\pi_{148}M$. The sphere action of $\beta_1$
commutes with $v_1^2$, so $v_1^3y=\beta_1v_1^2U=0$.
Exactness gives a lift $z\in\pi_{167}M_3$ with $\partial_3z=y$.

\section{Detection of the iterates}\label{sec:detect}

Fix the degree-$144$ self-map $V\colon\Sigma^{144}M_3\to M_3$
constructed in \cite[Corollary 4.7]{BS} by multiplication by
the class detected by $v_2^9$. This map induces $v_2^9$ on
$K(2)$-homology. For the lift $z\in\pi_{167}M_3$ just
constructed, set
\begin{equation}\label{eq:periodic-iterates}
  X_k=\partial_0\partial_3(V^k z)\in
          \pi_{153+144k}\Sph_{(3)},\qquad k\geq0 .
\end{equation}
Here $V^0$ is the identity, so $X_0=X$. Put $R=\Jconn$ and
let $T$ be $3$-completed $\mathrm{tmf}$. For $Y=R$ or $T$,
the notation $Y/(3,v_1^j)$
means $Y\wedge M_j$, where $\wedge$ is the smash product of
spectra. Let $h$ denote the unit $\Sph\to R$ and its induced
maps on homotopy groups, including those on cofibers obtained
by smashing with $M$ or $M_j$.

Write $T^\psi$ for the equalizer of $\psi^2$ and the identity
on $T$. The map $R\to T^\psi$ is an isomorphism on $\pi_d$
for $d\geq22$ \cite[Definition~2.2 and \S2.1]{CD}.
The operation $\psi^2$
fixes the torsion in $\pi_*T$
\cite[Theorem~C]{DavOps}\cite[Theorem~2.1]{CD}.
The permanence of $\xclass$ and the
detection criterion in \cite[Remark 3.14]{CD} give $h(X)\ne0$.
Write $\F\{g\}$ for the
one-dimensional $\F$-span of a generator $g$. The equalizer
sequence and $\pi_{153}T=0$ give
$\pi_{153}R=\F\{\partial(\beta\Delta^6)\}$. Thus
$h(X)\unitrel\partial(\beta\Delta^6)$. We also have
$\pi_{154}R\cong\F$
by $\pi_{155}T=0$ and $\pi_{154}T=\F\{\beta\Delta^6\}$.
These groups follow from the prime-three $\mathrm{tmf}$ calculation
of Hopkins and Mahowald, as presented by Bauer \cite[\S\S5--6]{Bau}.
See also the homotopy group calculation of Carrick, Davies, and
van Nigtevecht \cite[Corollary~C and Figure~A.1]{CDN}.

Choose the associative ring map $R/27\to R/9$ obtained by
smashing the compatible $\mathbb E_1$-algebra reduction
$\Sph/27\to\Sph/9$ with $R$ \cite[Theorem~1.5 with $n=1$]{BurMoore}.
The underlying unital map $\Sph/27\to\Sph/9$ is the usual
reduction, since $\pi_1(\Sph/9)=0$ at the prime $3$ and the
bottom-cell restriction determines the map.
As in \cite[\S3.3, before Lemma~3.22]{CD}, its fiber gives
$R/3$ the $R/27$-module structure used below, with
$R/27$-linear fiber inclusion $i\colon R/3\to R/27$.
Write $q_{27}\colon R\to R/27$ and $\partial_{27}\colon R/27\to\Sigma R$
for the mod-$27$ quotient and boundary. The factorization
$27=9\cdot3$ gives a map from the mod-$3$ to the mod-$27$
cofiber sequence via the identity on the source of multiplication
by $3$ and multiplication by $9$ on its target. The resulting
cofiber map is the fiber inclusion $i$, so
$iq_0=9q_{27}$ and $\partial_{27}i=\partial_0$.

The periodicity element $a=\widetilde{\partial\Delta^6}$ in
$\pi_{144}(R/27)$ is the unique class satisfying
$\partial_{27}(a)=\partial(\Delta^6)$
\cite[Proposition~3.17]{CD}.

The Moore boundaries give two comparisons, each up to sign.
Over $R$, applying $\partial_3$ carries the action of the
self-map $V$ to the action of $a$ on $R/3$. Over $T$, the action
of $a$ on the image of the $v_1$-boundary agrees with
multiplication by $\Delta^6$. Together these identities let us
compute the image of each $X_k$ in $R$ from the $\Delta^6$-action
on a lift of $\beta$ over $T$.

For $Y=R$ or $T$, let $a_N$ denote
its degree-$144$ action on $N=Y/3$. The action over $T$ is
obtained by extending scalars from $R$.

\begin{lemma}\label{lem:after-boundary}
Take $(Y,j)=(R,3)$ or $(T,1)$ and put $C=Y\wedge M_j$ and
$N=Y/3$. Let $U\colon\Sigma^{144}C\to C$ be $R\wedge V$ for
$(R,3)$ and multiplication by $\Delta^6$ for $(T,1)$. Then the
boundary $\delta=\partial_j\colon C\to\Sigma^{4j+1}N$ satisfies
\begin{equation}\label{eq:after-boundary}
  \delta U\unitrel(\Sigma^{4j+1}a_N)\,(\Sigma^{144}\delta).
\end{equation}
\end{lemma}

\begin{proof}
The low sphere stems give $\pi_1M=\pi_2M=0$ and
$\pi_0M\cong\mathbb Z/3$ by the mod-$3$ cofiber sequence.
Applying stable maps into $M$ to that sequence gives
$3\,\mathrm{id}_M=0$ and $[\Sigma M,M]=0$.
Since $3\,\mathrm{id}_{M_j}$ restricts to zero on $M$,
it factors through $\partial_j$, and hence $9\,\mathrm{id}_{M_j}=0$.

For $D=M$ or $M_j$, extend the identity on $D$ across the
bottom-cell inclusion $D\to(\Sph/9)\wedge D$ by a map
$\mu_D\colon(\Sph/9)\wedge D\to D$.
The difference between $\partial_j\mu_{M_j}$ and the composite
that applies $\partial_j$ first and then the suspended $\mu_M$
is zero on the bottom cell. It therefore factors through the
mod-$9$ boundary on $(\Sph/9)\wedge M_j$. Precomposition
by $\Sph/27\to\Sph/9$ multiplies that boundary by $3$, which kills
the difference since the target is a suspension of $M$.
Thus the resulting maps $(\Sph/27)\wedge D\to D$ commute
with $\partial_j$ up to homotopy. For $D=M$ the resulting map
agrees with the action on the fiber of $\Sph/27\to\Sph/9$.
Both restrict to the identity on $M$ and $[\Sigma M,M]=0$.

Smashing these maps with $Y$ gives $Y$-linear maps
$(Y/27)\otimes_Y(Y\wedge D)\to Y\wedge D$, where $\otimes_Y$
is the tensor product of $Y$-module spectra. Evaluate these maps
at the image of $a$ in $\pi_{144}(Y/27)$. Write
$A\colon\Sigma^{144}C\to C$ for the result on $C$. On $N$ the
result is $a_N$, so
$\delta A=(\Sigma^{4j+1}a_N)(\Sigma^{144}\delta)$.
Restricting the original map $(\Sph/27)\wedge M_j\to M_j$
along the bottom sphere of $M_j$ and then smashing with $Y$
gives a $Y$-linear map $Y/27\to C$ extending $Y\to C$.
This map agrees with the reduction followed
by the inclusion $Y/3\to C$, since $\pi_1C=0$. Indeed,
$\pi_1Y=0$ and $\pi_0Y$ is torsion-free for $Y=R,T$, and the
two cofiber sequences preserve vanishing in stem one.
The reduction calculation in \cite[Proposition~3.17 and its proof]{CD}
therefore gives the same bottom-cell value for $U$ and $A$ up to
a unit, using the $v_2^9$ class defining $V$ over $R$ and the
reduction of $\Delta^6$ over $T$. The value comes from
$\pi_{144}N$, which is killed by $3$, so this unit acts by a sign.

We next check $\pi_{145}N=\pi_{145}C=0$ and injectivity
of $v_1^j$ on $\pi_{144}N$. For $Y=T$, the homotopy calculation
\cite[Corollary~C and Figure~A.1]{CDN} and the mod-$3$ sequence
give $\pi_{133}(T/3)=\pi_{145}(T/3)=\pi_{146}(T/3)=0$.
For $n=132,140,144$, the group $\pi_nT$ is torsion-free and
$\pi_{n-1}T=0$. Let $c_4$ and $c_6$ denote the normalized
modular forms of weights four and six and their homotopy lifts.
In those three stems the reductions of the classes detected
by the following monomials form $\F$-bases of $\pi_n(T/3)$
\cite[Proposition~4.6]{HopkinsTMF}\cite[Corollary~C and Figure~A.1]{CDN}.
\begin{equation}\label{eq:tmf-reduction-bases}
  \begin{array}{ccl}
    n & \text{detecting monomials} & \text{range}\\
    132 & c_4^{15-3b}c_6\Delta^b & 0\leq b\leq5\\
    140 & c_4^{16-3b}c_6\Delta^b & 0\leq b\leq5\\
    144 & c_4^{18-3b}\Delta^b & 0\leq b\leq6
  \end{array}
\end{equation}
By $T$-linearity and \cite[Equation~(2.7)]{CD}, on reductions of
classes from $T$ the maps $v_1^2$ and $v_1^3$ act by $c_4$ and
$c_6$, respectively. The products are independent after
reduction of the modular forms, using $c_6^2=c_4^3$ modulo
three, so the products are independent in $\pi_*(T/3)$.
This proves injectivity of $v_1^j$ on $\pi_{144-4j}(T/3)$
and on $\pi_{144}(T/3)$ for $j=1,3$.

The first injection and $\pi_{145}(T/3)=0$ give
$\pi_{145}(T\wedge M_j)=0$. Also $\pi_{146}(T\wedge M_3)=0$
by the vanishing in stems $133$ and $146$ of $T/3$.
The high-stem comparison $R\to T^\psi$ applies after
smashing with $M$ or $M_3$ by their defining cofiber sequences.
The equalizer sequence after these smashes gives
$\pi_{145}(R/3)=\pi_{145}(R\wedge M_3)=0$ and an injection
$\pi_{144}(R/3)\to\pi_{144}(T/3)$. Naturality now gives the
required injectivity of $v_1^3$ in stem $144$ of $R/3$.

Choose $\epsilon=\pm1$ so that $U-\epsilon A$ restricts to
zero on the bottom cell. Since $\pi_{145}C=0$, it also restricts
to zero along $\Sigma^{144}N\to\Sigma^{144}C$. It follows
that $U-\epsilon A$ factors as $\Sigma^{144}\delta$ followed by a
$Y$-linear map $e\colon\Sigma^{145+4j}N\to C$.
The desuspension of $\delta e$ is a map $\Sigma^{144}N\to N$
whose bottom-cell value lies in the kernel of $v_1^j$ on
$\pi_{144}N$ by exactness. That kernel is zero. Since
$\pi_{145}N=0$, a $Y$-linear map $\Sigma^{144}N\to N$
is determined by its bottom-cell value, so $\delta e=0$.
This gives $\delta U=\epsilon\delta A$ as claimed.
\end{proof}

Let $b_R=h(y)\in\pi_{154}(R/3)$ and write $a^k b_R$ for
successive applications of $a_N$ to $b_R$. Iterating
\eqref{eq:after-boundary} gives
$h(X_k)\unitrel\partial_0(a^k b_R)$.

Choose a lift $\widetilde\beta\in\pi_{11}(T/3)$ with
$\partial_0(\widetilde\beta)=\beta$.
Bauer's calculation \cite[Figure~1]{Bau} gives
$\pi_{14}T=\pi_{15}T=0$, so $\pi_{15}(T/3)=0$
\cite[Figure~1]{CD}. Thus
$v_1\widetilde\beta=0$, and exactness gives a lift of
$\widetilde\beta$ along the $v_1$ boundary from $\pi_{16}T/(3,v_1)$.
The $j=1$ case of \eqref{eq:after-boundary} identifies the
$a_N$ action on the image of the boundary with multiplication
by $\Delta^6$ up to sign. Both preserve the image, so iteration gives
$a^k\widetilde\beta\unitrel\Delta^{6k}\widetilde\beta$
in $\pi_*(T/3)$.

The equalizer boundary $\Sigma^{-1}T\to T^\psi$ is $R$-linear.
Tensoring over $R$ with the $R/27$-module $R/3$ gives the
$R/27$-linear map $\Sigma^{-1}(T/3)\to T^\psi/3$ used in
\cite[proof of Lemma~3.23]{CD}. The map $R/3\to T^\psi/3$
is an isomorphism on homotopy in stems at least $23$, by the
mod-$3$ cofiber sequences and the corresponding comparison for
$R$. Write $\bar\partial$ for the induced boundary on homotopy
in this range, transported to $\pi_*(R/3)$ by this isomorphism.
It commutes with the action of $a$ there.
By naturality of the mod-$3$ cofiber sequences, the lift
$b_0=\bar\partial(\Delta^6\widetilde\beta)$ of
$\partial(\beta\Delta^6)$ in $\pi_{154}(R/3)$ satisfies
\begin{equation}\label{eq:periodic-boundary}
\partial_0(a^k b_0)\unitrel
 \partial\bigl(\beta\Delta^{6(k+1)}\bigr).
\end{equation}
Here the boundary $\partial_0$ on $T/3$ is linear over $T$.

After rescaling $b_0$ by a unit, exactness gives
$b_R-b_0=q_0u_R$ for a class $u_R\in\pi_{154}R$.
Since $3u_R=0$, the identity $iq_0=9q_{27}$ gives $iq_0u_R=0$.
Using the $R/27$-linearity of $i$ and $\partial_0=\partial_{27}i$, we get
$\partial_0(a^k(b_R-b_0))=\partial_{27}(a^k i q_0u_R)=0$.
This proves
\begin{equation}\label{eq:periodic-detection}
h(X_k)\unitrel\partial\bigl(\beta\Delta^{6(k+1)}\bigr)\ne0
\end{equation}
for all $k$. The final nonvanishing follows because
$\beta\Delta^{6(k+1)}$ spans the order-three group in stem
$154+144k$ of $T$ \cite[Corollary~C and Figure~A.1]{CDN},
on which $\psi^2-1$ is zero \cite[Theorem~2.1]{CD}.

Let $\nu T$ denote the $\BP$-synthetic analogue of $T$
\cite[Definition~4.3]{Pst}\cite[Example~9.2 and Construction~9.3]{BHS}.
To determine the filtration, use $(\nu T)^\psi$, the equalizer
of $\nu(\psi^2)$ and the identity
on $\nu T$ \cite[Definition~2.10 and Equation~(2.11)]{CD}.
Its filtration comes from the tower of multiplication by $\tau$
\cite[\S2.2, Equation~(2.13)]{CD}.
Inverting $\tau$ recovers $T^\psi$ and its homotopy boundary,
since $\tau$-inversion is exact and takes $\nu T$ to $T$
\cite[Theorem~9.12(3)]{BHS}\cite[proof of Proposition~2.12]{CD}.
The equalizer boundary decreases the stem by one and raises
filtration by one. The synthetic model of $R$ maps to
$(\nu T)^\psi$ and induces isomorphisms on bigraded homotopy
in stems at least $22$, in every filtration
\cite[Definition~2.10, Equation~(2.11), and Proposition~2.12]{CD}.
Each filtration subgroup is the image of the corresponding
stage of the tower in homotopy. Thus this comparison identifies
the filtrations on $\pi_*R$ and $\pi_*T^\psi$ in those stems.

The bigraded homotopy of $\nu T$ modulo $\tau$ gives the
Adams--Novikov $E_2$-page of $T$ \cite[Corollary~9.13]{BHS}.
On that page,
$\mathcal E_2(154+144k,2;T)$ is spanned over $\F$ by
$\beta\Delta^{6(k+1)}$, which is permanent as a product of
$\beta$ and powers of $\Delta^3$
\cite[Corollary~D and Figure~A.1]{CDN}. The comparison from the
$3$-local $\BP$-Adams tower of $\mathrm{tmf}$ to the tower of
$T$ is degreewise $3$-completion. Indeed, each smash power of
the fiber of $\Sph_{(3)}\to\BP$ is connective and of finite
type, so for fixed homotopy degree one can replace it by a
finite skeleton before smashing with $\mathrm{tmf}$ or $T$.
The tower groups are degreewise finitely generated before
completion. Tensoring them with the $3$-adic integers is exact,
so the order-three group on $E_2$ is unchanged.
The whole group survives to $E_\infty$ and detects the order-three
group of $\pi_{154+144k}T$. Since $\psi^2$ is the identity on
that homotopy group \cite[Theorem~2.1]{CD}, it is already the
identity on this $E_2$ group.
As in Section~\ref{sec:proof}, the mod-$\tau$ exact sequence
injects the cokernel of $\psi^2-1$ via the equalizer boundary.
This gives a nonzero class in bidegree $(153+144k,3)$
in the spectral sequence for $T^\psi$. Naturality makes it a
cycle on every page. The $E_2$-page for $T$ has nonnegative
filtration and internal degrees divisible by four
\cite[Corollary~D]{CDN}. Both properties pass to the equalizer
page, since the boundary shifts the bidegree by $(-1,1)$.
Hence no differential can hit filtration three.
The surviving class represents the image of $h(X_k)$ in
$\pi_*T^\psi$ by \eqref{eq:periodic-detection} and compatibility
of the spectral sequence boundary with the homotopy boundary.
The comparison with $T^\psi$ therefore places $h(X_k)$ in exact
filtration three on $R$. The unit from the sphere to $R$ is
filtered \cite[\S2.2, after Equation~(2.13)]{CD}, so $X_k$ has filtration
at most three. Since $153+144k\equiv1\pmod4$, the first possible
filtration in the sphere in that stem is three. Thus $X_k$ has
exact filtration three. It has exact order three because it
is nonzero and lies in the image of $\partial_0$. This completes
the proof of Theorem~\ref{thm:main}.

\section*{Acknowledgment}
The author is supported by the NSFC Grant No.~12405001. The author thanks AI for its help.

\end{document}